\documentclass[12pt,a4paper]{amsart}
\usepackage{graphicx}
\usepackage{amssymb}
\usepackage{setspace}
\newtheorem{thm}{Theorem}[section]
\newtheorem{fact}[thm]{Fact}

\newtheorem{lem}[thm]{Lemma}
\newtheorem{prop}[thm]{Proposition}
\theoremstyle{definition}
\newtheorem{defn}[thm]{Definition}
\theoremstyle{remark}

\numberwithin{equation}{section}

\newcommand{\<}{\langle}
\renewcommand{\>}{\rangle}
\begin{document}

\title[]{One variable equations  in  torsion-free hyperbolic groups}%
\thanks{This research was supported by the "ANR" Grant  in the framework of the project "GGL"(ANR-05-JC05-47038).}%

\author{Abderezak OULD HOUCINE}%
\address{Universit\'e de Lyon;
Universit\'e Lyon 1;
INSA de Lyon, F-69621;
Ecole Centrale de Lyon;
CNRS, UMR5208, Institut Camille Jordan,
43 blvd du 11 novembre 1918,
F-69622 Villeurbanne-Cedex, France}%
\email{ould@math.univ-lyon1.fr}%
\email{}%

\thanks{}%
\subjclass{}%
\keywords{}%


\newcommand{\cal}[1]{\mathcal{#1}}
\newcommand{\di}{{\rm{dim}}}
\newcommand{\Di}{{\rm{Dim}}}
\newcommand{\Rc}{{\rm{RC}}}
\newcommand{\Rm}{{\rm{RM}}}
\newcommand{\tp}{{\rm tp}}

\begin{abstract}  Let $\Gamma$ be a torsion-free hyperbolic group. We show that     the set of solutions of any  system of equations with one variable  in $\Gamma$  is a finite union of points and cosets of centralizers if and only if any two-generator subgroup of $\Gamma$ is free. \end{abstract}
\maketitle

\section{Introduction}

Equations with one variable in free groups have been studied by Lyndon \cite{Lyndon-equa}, Lorents \cite{Lorents1, Lorents2},   and Appel \cite{Appel-equation},  among others, and the  conclusion is that  the set of solutions of a finite system of equations  with one variable  is a finite union of points and cosets of centralizers.  However,  Lorents announced his result without proof and the proof of Appel contains a gap \cite{baumslag-review}.  In \cite{chis-rem}, Chiswell and Remeslennikov gave a   proof of this result, by using coordinate groups and Lyndon length functions in ultrapowers of free groups. In this paper we shall be concerned with  a description of equations with one variable in a more  larger class  of  groups.

\begin{thm} \label{thmmain}  Let $\Gamma$ be a nonabelian torsion-free hyperbolic group such that any two-generator subgroup of $\Gamma$ is free. Then   the set of solutions of a  system of equations with one variable  in $\Gamma$  is a finite union of points and cosets of centralizers. \end{thm}

We notice that nonfree torsion-free hyperbolic groups whose two-generator subgroups are free exist. For instance by taking a nonfree hyperbolic group which is a limit group of free groups, we obtain such examples.  The precedent  theorem is no longer true if we drop the assumption that two-generator subgroups of $\Gamma$ are free,  and in fact we have the following equivalence.

\begin{thm} \label{thmmain2}  Let $\Gamma$ be a nonabelian torsion-free hyperbolic group. Then the following properties are equivalent:

$(1)$  the set of solutions of a  system of equations with one variable  in $\Gamma$  is a finite union of points and cosets of centralizers; 

$(2)$ any two-generator subgroup of $\Gamma$ is free. 

\end{thm}

 As a consequence of Theorem \ref{thmmain} is that any quantifier-free formula in a torsion-free hyperbolic group satisfying the hypothesis of the theorem, is a boolean combination of cosets of centralizers. It follows in particular that any proper subgroup of a torsion-free hyperbolic group,  under  the hypothesis of the theorem, which is definable by a quantifier-free formula,  is cyclic. However this property is steal true in any torsion-free hyperbolic group and has a simple proof (see the end of the appendix). 
 
Our approach to prove  Theorem \ref{thmmain} is to use coordinate groups of varieties as in \cite{chis-rem},  and the structure of restricted $\Gamma$-limit groups  obtained from  Sela's work on limit groups of torison-free hyperbolic groups \cite{Sela-hyp}.  In the next section we prove the main result, while the appendix is devoted to the proof  of an intermediate result on the structure of  restricted $\Gamma$-limit groups.

\section{Equations with one variable}

Let $G$ be a  group. A
\emph{$G$-group} $H$ is a group having an isomorphic fixed copy of
$G$ which we will identify with $G$. A homomorphism $h : H_1
\rightarrow H_2$ between two $G$-groups is called a
\emph{$G$-homomorphism} if for any $g \in G$, $h(g)=g$.  A \textit{$G$-isomorphism}  is  defined analogously and  we use the notation $H_1 \cong_G H_2$.

If $X$ is
set, we denote by $G[X]$ the free product $G*F(X)$ where $F(X)$ is
the free group on $X$. If  $X=\bar x=\{x_1, \dots, x_n\}$ we use the notation  $G[\bar
x]$.

For an element $w(\bar x) \in G[\bar x]$ and a tuple  $\bar g
=(g_1, \dots, g_n)\in G^n$ we denote by $w(\bar h)$ the element of
$G$ obtained by replacing each $x_i$ by $g_i$ ($1 \leq i \leq n$). A \textit{variety} in $G^n$ is a set of the form
$$
V(S)=\{\bar g \in G^n ~|~w(\bar g)=1 \hbox{ in $G$ for all }w \in S\},
$$ 
for some $S \subseteq G[ \bar x]$. For any  $S \subseteq G[ \bar x]$ we use the notation $S(\bar x)=1$ as an abreviation for the system of equations $\{w(\bar x)=1| w\in S\}$. The group $G$ is called \emph{equationally noetherian} if for any $n \geq 1$ and
any  subset $S$ of $G[\bar x]$ there exists a finite subset $S_0
\subseteq S$ such that ${V(S)=V(S_0)}$. A subset of $G^n$ is
\textit{closed} if it is the intersection  of
finite union of  varieties. This defines a topology on $G^n$, called the \textit{Zariski topology}.  Then $G$ is equationally noetherian if and only if for any $n \geq 1$, the Zariski topology on $G^n$ is noetherian \cite{Baum-Rem-algebraic}.  

If $G$ is equationally noetherian, then for any variety $V(S)$ in $G^n$ one associates to it its irreducibles components, which are also varieties. For more details on these notions we refer the reader  to \cite{Baum-Rem-algebraic}.

\cite[Theorem 1.22]{Sela-hyp} states that any system of equations (without parameters) in finitely many variables is equivalent in a trosion-free hyperbolic group to a finite subsystem.  This property is equivalent, when the group under consideration $G$ is finitely generated, to the fact that $G$ is equationally noetherian. Hence, we have the following.

\begin{fact} \cite[Theorem 1.22]{Sela-hyp} A torsion-free hyperbolic group is equationally noetherian. \qed
\end{fact}

It follows that the set of solutions of a system of equations is a finite union of irreducible varieties, and thus the study of such sets is reduced to the study of irreducible varieties.  This section is devoted to prove the following theorem which is a generalization of \cite[Theorem 5.3]{chis-rem}.

\begin{thm} \label{thm1} Let $\Gamma$ be a nonabelian torsion-free hyperbolic group such that any two-generator subgroup of $\Gamma$  is free. If $V$ is a proper nonempty irreducible variety in  $\Gamma$ then either $V$ is  a singleton  or $V$ is a coset of a centralizer. \end{thm}

Theorem \ref{thmmain} is a mere consequence of  Theorem \ref{thm1}.  Concerning quantifier-free formulas, since $\Gamma$ is a nonabelian CSA-group, there exist two elements $c, d \in \Gamma$ such that $C_\Gamma(c,d)=C_\Gamma(c) \cap C_\Gamma(d)=1$,  and thus  we can write $V=aC_\Gamma(c) \cap aC_\Gamma(d)$ when $V=\{a\}$.

Let $G$ be a group and $S$ a subset of $G[x_1, \dots, x_n]$. We let 
$$
G_S(\bar x)=\<G[x_1,\dots, x_n]| w(x_1, \dots, x_n)=1, \;\; w\in S\>,
$$
and we let 
$$
\bar S=\{w(\bar x) \in G[\bar x] | G \models \forall \bar x (S(\bar x)=1 \Rightarrow w(\bar w)=1)\}.
$$ 

The group $G_{\bar S}(\bar x)$ is called the \textit{coordiante group} associated to $S$ or to $V(S)$. We notice that for any $w \in G[\bar x ]$, $G_{\bar S}(\bar x) \models w(\bar x)=1$ if and only if $w \in \bar S$.   

In order to prove  the above theorem we shall need Lemma \ref{lem-reduction} below,  which connects  the structure of a variety to the structure of its coordiante group.  First,  we prove the following  technical proposition of independent interest. 







\begin{prop} \label{lem-hnn}  Let $G$ be a group, $H$ a subgroup of $G$ and suppose that  $G$ is generated by $H \cup\{s\}$ for some $s \in G$. 

$(1)$ Let  $G=\<H, t| [A,t]=1\>$ and suppose that: 

 $(i)$ $A$ is malnormal in $H$;

 $(ii)$ for any $u \in H$ either $u \in A$ or $\<u,A\>$ is the natural free product $\<u\>*A$. 

Then  there exist $h_1, h_2 \in H$ such that  $s^{\pm 1}=h_1th_2$. 

$(2)$ Let  $G=\<H, t| A^t=B\>$ and suppose that: 

 $(i)$ $A$ and $B$ are  malnormal in $H$ and $G$ is seperated; 

 $(ii)$ for any $u \in H$ either $u \in A$ or $\<u,A\>$ is the natural free product $\<u\>*A$ and similarly for $B$.  

Then  there exist $h_1, h_2 \in H$ such that  $s^{\pm 1}=h_1th_2$. 

$(3)$ Let  $G=H*_{A=B}K$ and suppose that: 

 $(i)$ $A$ is malnormal in $H$ and  $B$ is malnormal in $K$; 

 $(ii)$ for any $u \in H$ either $u \in A$ or $\<u,A\>$ is the natural free product $\<u\>*A$, and similarly for $B$ in  $K$.  

Then  there exist $h \in H, k \in K$ such that  $s^{\pm 1}=hk$. 
\end{prop}

\proof $\;$

$(1)$ If $v(\bar x,y)$ is a word in the free group with basis $\bar x \cup\{y\}$,  we denote by $exp_y(v)$ the exponent sum of $y$ in $v$. 

Let $s'$ be a cyclically reduced conjugate of $s$. Since $H \cup\{s\}$ generates $G$, there is a word $w(\bar x, y)$ such that $t=w(\bar h, s)$ for some tuple $\bar h$ of $H$.  By the abelianization of $G$, we have $exp_t(w(\bar h, s))=exp_t(s')exp_s(w)$. Therefore $exp_t(s')=\pm 1$.  Hence the number of occurrences of $t$ in $s'$ is odd. We have  $s=s'^g$ for some $g \in G$.  We claim now that $g \in H$ and $|s'|=1$,  where $|.|$ denotes the length of normal forms. 

Suppose towards a contradiction that $|s'|  \geq 1$ and thus $|s'|\geq 3$. Using the fact that $exp_t(s')=\pm 1$ and the malnormality of $A$, a simple count shows that 
$$
|s^2|>|s|, \leqno (1)
$$
and using also calculations with normal forms,  we get  for any $h,h' \in H$, with $h \neq 1, h'\neq 1$, that $$|hs^{\pm 1}h's^{\pm 1}|> |hs^{\pm 1}|, |h's^{\pm 1}|. \leqno (2)$$  

Using (1) and (2),  by \cite[Lemma 4.2]{AHoucine}, we get that or any sequence $h_1, \dots, h_n$ of nontrivial elements of $H$, for any sequence $\varepsilon_0, \dots, \varepsilon_{n}$ of $\mathbb Z$, $\varepsilon_i \neq 0$,  
$$
|s^{\varepsilon_0}h_1 s^{\varepsilon_1}h_2 \cdots h_ns^{\varepsilon_{n}}|\geq |s|>1,
$$
and thus $t \not \in \<H,s\>$;  which is a contradiction. 

Therefore $|s'|=1$ and we write $s'=ut^{\varepsilon}v$, where $\varepsilon=\pm 1$. To simplify, we may assume that $\varepsilon=1$.  Write $g = h_0t^{\varepsilon_1}h_1\cdots t^{\varepsilon_n}h_n$ in normal form.  Replacing $s$ by $h_n^{-1}sh_n$ and $s'$ by $h_1^{-1}utvh_1$ we may assume without loss of generality that $h_0=h_n=1$. 

We claim now that $h_1, \dots, h_{n-1} \in A$. Suppose that for some $i$,   $h_i \in A$. Then proceeding as above, a simple count with normal forms, shows that  for any $h,h' \in H$, with $h \neq 1, h'\neq 1$, that $|hs^{\pm 1}h's^{\pm 1}|> |hs^{\pm 1}|, |h's^{\pm 1}|$, and we get  a contradiction as above  by \cite[Lemma 4.2]{AHoucine}. 

Hence we get $g=a t^p$ for some $p \in \mathbb Z$ and $a \in A$. Replacing again $u$ be $a^{-1}u$ and $v$ by $va$, we may assume that $a=1$.  Hence $s=t^{-p}utvt^p$. 

We claim that either $v \in A$ or $u \in A$.  Suppose that $v \not \in A$ and $u \not \in A$. Then proceeding as above, we see also that a simple calculation with normal forms, shows that  for any $h,h' \in H$, with $h \neq 1, h'\neq 1$, that $|hs^{\pm 1}h's^{\pm 1}|> |hs^{\pm 1}|, |h's^{\pm 1}|$, which is a contradiction by \cite[Lemma 4.2]{AHoucine}. 

Hence $v \in A$ or $u \in A$. We treat only the case $v \in A$,  the other case being similar.   Replacing again $u$ by $ua$, we may assume that $v=1$. Therefore $s= t^{-p}ut^{p+1}$.  Clearly by $(ii)$ we have $\<u,A\>=\<u\>*A$. 

We claim that for any sequence $h_1, \dots, h_n$ of nontrivial elements of $H$, for any sequence $\varepsilon_0, \dots, \varepsilon_{n}$ of $\mathbb Z$, $\varepsilon_i \neq 0$,   the normal form of the product 
$$
s^{\varepsilon_0}h_1 s^{\varepsilon_1}h_2 \cdots h_ns^{\varepsilon_{n}},
$$
is of the form 
$$
t^{\delta_1}d_1 \cdots t^{\delta_p}d_pt^{q},
$$
where $\delta_i=\pm 1$,  $q \in \{p, p+1\}$,  $d_i \in H$,  and $d_p \in \<u,A\>$ with the property that the last element of the normal form of $d_p$, with respect to the structure   $\<u,A\>=\<u\>*A$,  is $u^{\pm 1}$. 

The proof is by induction on $n$ and the detailled verification is left to the reader. 

Hence we conclude that for any sequence $h_1, \dots, h_n$ of nontrivial elements of $H$, for any sequence $\varepsilon_0, \dots, \varepsilon_{n}$ of $\mathbb Z$, $\varepsilon_i \neq 0$,  
$$
|s^{\varepsilon_0}h_1 s^{\varepsilon_2}h_2 \cdots h_ns^{\varepsilon_{n}}|\geq 2,
$$
and thus $t \not \in \<H,s\>$; a final contradiction.

(2) This case is similar to (1).  Proceeding as above, we conclude that $s=g^{-1}utvg$, and we suppose that $g \not \in H$. Then, as before we may assume that  $g=t^{\varepsilon_1}h_1 \cdots h_{n_1}t^{\varepsilon_n}$. At this stage, by using the fact that $G$ is seperated, we get $g=t^{\pm 1}$ and we assume without loss of generality that $g=t$. Hence $s=t^{-1}utvt$.  Then, as above, we may assume that $v \in B$ and, without loss of generality $v=1$ and thus  $s=t^{-1}ut^2$.  Then as before, for any sequence $h_1, \dots, h_n$ of nontrivial elements of $H$, for any sequence $\varepsilon_0, \dots, \varepsilon_{n}$ of $\mathbb Z$, $\varepsilon_i \neq 0$,   the normal form of the product 
$$
s^{\varepsilon_0}h_1 s^{\varepsilon_2}h_2 \cdots h_ns^{\varepsilon_{n}},
$$
is of the form 
$$
t^{\delta_1}d_1 \cdots t^{\delta_p}d_pt^{q},
$$
where $\delta_i=\pm 1$,  $q \geq 1$, $d_i \in H$,  and $d_p \in \<u,A\>$ with the property that the last element of the normal form of $d_p$, with respect to the structure   $\<u,A\>=\<u\>*A$,  is $u^{\pm 1}$. The proof is by induction on $n$ and the detailled verification is left to the raider. We conclude that $t \not \in \<H,s\>$; a final contradiction.

(3) This case is also similar to (1) and (2).    Write  $s =y_1y_2 \dots y_n$ in normal form.  We claim that $n \leq 2$.

Suppose first by contradiction that $n \geq 4$.   Using calculations with normal forms, we find that $| h y^{\pm 1} h ' y^{\pm 1}|>|h y^{\pm 1}|, |h' y^{\pm 1}|$ for any $h,h'$ with $h \neq 1$ and $h' \neq 1$. Hence by \cite[Lemma 4.2]{AHoucine},  for any nontrivial elements $h_1, \dots, h_n \in H$,  for any sequence $p_1, \dots, p_n$ of $\mathbb Z$, $p_i \neq 0$, we have 
$$
|h_1 y ^{p_1} \cdots h_n y^{p_n}|>|y|>1,
$$
which is clearly a contradiction.  Therefore $n \leq 3$.  Suppose that $n=3$.  We treat only the case $y_1 \in H$, the other case being similar. 

We claim now that for any sequence $h_1, \dots, h_n$ of nontrivial elements of $H$, for any sequence $\varepsilon_0, \dots, \varepsilon_{n}$ of $\mathbb Z$, $\varepsilon_i \neq 0$,   the normal form of the product 
$$
s^{\varepsilon_0}h_1 s^{\varepsilon_2}h_2 \cdots h_ns^{\varepsilon_{n}},
$$
is of the form 
$$
d_1 \cdots d_p d_{p+1},
$$
where $ p \geq 2$,  $d_{p+1} \in \{y_1^{\pm1}, y_3^{\pm 1}\}$,  and $d_p \in \<y_2,B\>$ with the property that the last element of the normal form of $d_p$, with respect to the structure   $\<y_2,B\>=\<y_2\>*B$,  is $y_2^{\pm 1}$. Which is a contradiction.   \qed

\begin{lem} \label{lem-reduction} Let $G$ be a group and $S \subseteq G[x]$. 

\smallskip
$(1)$ If $G_{\bar S}(x) \cong_G G$ then $V(S)$ is a singleton. 

$(2)$ If $G_{\bar S}(x) \cong_G G*\mathbb Z$ then $V(S)=G$. 

$(3)$ If $G_{\bar S}(x) \cong_G \<G,t | [A,t]=1\>$, where $A$ is a nontrivial malnormal cyclic subgroup of $G$,  and $G$ satisfies the property that any two-generator subgroup of $G$ is free, then $V(S)=uC_G(A)^v$ for some $u,v \in G$. 
\end{lem}

\proof

$(1)$ Let $h : G_{\bar S}(x) \rightarrow G$ be a $G$-isomorphism. Then $h(x) \in G$ and thus $x \in G$. Hence $xg^{-1} \in \bar S$ for some $g \in G$. Therefore $V(S)=\{g\}$.

$(2)$  Let $h : G_{\bar S}(x) \rightarrow G*\mathbb Z$ be a $G$-isomorphism. Clearly $h(x) \not \in G$. Hence the subgroup $\<G, h(x)\>$ is the natural free product $ G*\<h(x)\>$.  Therefore if $G_{\bar S}(x) \models w(x)=1$, where $w \in S$, then $G*\<t|\> \models w(t)=1$. Hence $G \models \forall t w(t)=1$ and thus $V(S)=G$. 

$(3)$  Let $h : G_{\bar S}(x) \rightarrow \<G, t |[A,t]=1\>$ be a $G$-isomorphism. Since $G \cup\{h(x)\}$ generates the HNN-extension under consideration,  $h(x)=u_0t^\varepsilon v$ where $ \varepsilon =\pm 1$ and $u_0,v \in G$ by Proposition  \ref{lem-hnn}.  Without loss of generality we may assume that $\varepsilon =1$.

We claim that  $V(S)=uC_G(A)^v$, where $u=u_0v$. Let $g \in G$ be a solution of the system $S(x)=1$. Then there exists a $G$-homomorphism $f : G_{\bar S}(x) \rightarrow G$ such that $f(x)=g$.  Hence $g=u_0f(t)v$. Since $[t,a]=1$ for all $a \in A$ we get $[f(t),a]=1$ for all $a \in A$ and thus $g \in uC_G(A)^v$.

We have $G_{\bar S}(x)=\<G, x | x=u_0tv, [A,t]=1\>$. Hence for any $w \in S$, $w(x)=1$ is a consequence of the precedent presentation.  Let $g \in uC_G(A)^v$. Then, using the precedent presentation, there exists a $G$-homomorphism $f : G_{\bar S}(x) \rightarrow G$ such that $f(x)=g$.  Hence for any $w \in S$, $w(g)=1$, by the precedent observation. We conclude that $V(S)=uC_G(A)^v$ and this terminates the proof.  \qed

Using Lemma \ref{lem-reduction}, the proof of Theorem \ref{thm1} is reduced  to the proof of  the following theorem which is a generalization of \cite[Theorem 5.1]{chis-rem}.

\begin{thm} \label{thm2}Let $\Gamma$ be a nonabelian torsion-free hyperbolic group such that any two-generator subgroup of $\Gamma$ is free. The coordinate group $\Gamma_{\bar S}(x)$  of the  nonempty  irreducible  variety $V(S)\subseteq \Gamma$
satisfies one of the following:

\smallskip
$(1)$  $\Gamma_{\bar S}(x) \cong_\Gamma \Gamma$; 

$(2)$  $\Gamma_{\bar S}(x) \cong_\Gamma \Gamma*\mathbb Z$;  

$(3)$  $\Gamma_{\bar S}(x) \cong_\Gamma \<\Gamma, t | [u, t] = 1\>$,  for some nontrivial element $u$ in  $\Gamma$. 
\end{thm}

The remainder of this section is devoted to prove Theorem \ref{thm2}.  In the sequel we let $\Gamma$ to be a fixed nonabelian torsion-free hyperbolic group.

\begin{defn} 
$\;$

$\bullet$ A sequence of homomorphisms $(f_n)_{n \in \mathbb N}$ from  $H$ to $\Gamma$ is called \textit{stable} if for any $h \in H$ either $f_n(h)=1$ for all but finitely many  $n$, or   $f_n(h) \neq 1$ for all but finitely many $n$.  The \textit{stable 
kernel} of $(f_n)_{n \in \mathbb N}$, denoted $Ker_{\infty}((f_n)_{n \in \mathbb N})$,  is the set of elements $h \in H$ such that $f_n(h)=1$ for all but finitely many $n$. 

\smallskip
$\bullet$  A \textit{restricted $\Gamma$-limit group} is a $\Gamma$-group $G$ such that  there exists a $\Gamma$-group  $H$  and a stable sequence of $\Gamma$-homomorphisms $(f_n)_{n \in \mathbb N}$ from $H$ to $\Gamma$ such that $G=H/ Ker_{\infty}((f_n)_{n \in \mathbb N})$. 
\end{defn}

The proof of the following lemma is straightforward and it is left to the reader.

\begin{lem}\label{lem1}  Let $S \subseteq \Gamma[x]$ such that $V(S)$ is irreducible and nonempty.

$(1)$   The group   $\Gamma_{\bar S}(x)$ is a $\Gamma$-group and for any finite subset $A$ of $\Gamma_{\bar S}(x)$ such that $1 \not \in A$ there exists a $\Gamma$-homomorphism $f : \Gamma_{\bar S}(x) \rightarrow \Gamma$ such that $1 \not \in f(A)$.

$(2)$  Let $\varphi : \Gamma_{\bar S}(x) \rightarrow L$ be a $\Gamma$-epimorphism where $L$ is a restricted $\Gamma$-limit group. Then  there exists  $U \subseteq \Gamma[x]$ such that $V(U)$ is irreducible and nonmepty and $L=\Gamma_{\bar U}(x)$.  \qed
\end{lem}

It follows in particular, by Lemma \ref{lem1}(1),  that if $V(S)$ is irreducible and nonempty then   $\Gamma_{\bar S}(x)$ is a restricted $\Gamma$-limit group.   Lemma \ref{lem1}(1) implies also that $\Gamma$ is existentially closed  in $\Gamma_{\bar S}(x)$.


\begin{defn} \cite[Definition 4.16]{Champ-Guirardel} Let $G$ be a group which is the fundamental group of a graph of groups $\Lambda$. Let $H$ be a nontrivial elliptic subgroup of $G$ with respect to $\Lambda$. The \textit{elliptic abelian neighbourhood} of $H$ is the subgroup $\hat H$ generated by all the elliptic elements of $G$ which commute with a nontrivial element of $H$.  \qed
\end{defn}


\begin{defn}  \label{weakly-constructible} A restricted $\Gamma$-limit group $G$ is said \textit{weakly constructible} if one of the following cases holds: 

$(1)$  $G=H*_CK$, where $\Gamma \leq H$ and $C$ is a notrivial cyclic group and $K$ is noncyclic; 

$(2)$ $G=\<H, t| C_1^t=C_2\>$, $\Gamma \leq H$ and $C_1$ is a notrivial cyclic group,  and there exists  a proper quotient restricted $\Gamma$-limit group $L$ of $G$ where the corresponding  $\Gamma$-epimorphism $\varphi : G \rightarrow L$  is one-to-one in restriction to the elliptic abelian neighbourhood of $H$. 
\end{defn}

We will use the following theorem, which is sufficient for our purpose, and   whose proof proceeds in a similar way to that  of  \cite{Sela-Diophan1,Champ-Guirardel}. For completeness,   the proof is given in the appendix.

\begin{thm} \label{thm-princip2} Let $G$ be a restricted $\Gamma$-limit group. If $G$ is not $\Gamma$-isomorphic to $\Gamma$ and if it is freely indecomposable relative to $\Gamma$ then $G$   is weakly constructible. \qed
\end{thm} 

\smallskip
\noindent
\textbf{Proof of Theorem \ref{thm2}.}  

Let $S \subseteq \Gamma[x]$ such that $V(S)$ is irreducible and nonempty. We may assume that   $\Gamma_{\bar S}(x)$ is not $\Gamma$-isomorphic to $\Gamma$ and it is freely indecomposable relative to $\Gamma$. By Lemma \ref{lem1}, every proper quotient of $\Gamma_{\bar S}(x)$, which is a restricted $\Gamma$-limit group,  is of the form $\Gamma_{\bar U}(x)$. Hence, by the descending chain condition on $\Gamma$-limit groups, we may assume that the theorem holds for all proper quotients of $\Gamma_{\bar S}(x)$ which are  restricted $\Gamma$-limit groups. By Theorem \ref{thm-princip2}, we treat the  two cases (1) and (2) of  Definition \ref{weakly-constructible}.

Since any two-generator subgroup of $\Gamma$ is free, we have the following claim whose proof proceeds in a similar way to that of \cite[Claim 4.25]{Champ-Guirardel} and it is left to the reader. 

\smallskip
\noindent Claim 1. \textit{Let $G$ be a $\Gamma$-limit group and $T_1, T_2$ two abelian subgroups of $G$. Then either $\<T_1,T_2\>$ is abelian or $\<T_1, T_2\>=T_1*T_2$.} \qed

\smallskip

We now prove  the following claim.

\smallskip
\noindent 
Claim 2. \textit{Suppose that $\Gamma_{\bar V}(x)=\<H, t|C_1^t=C_2\>$, where $\Gamma \leq H$ and $C_1$ is cyclic,  is a restricted $\Gamma$-limit group which  is freely indecomposable relative to $\Gamma$ and such that: }

$(i)$\textit{ there exists  a proper quotient restricted $\Gamma$-limit group $L$ of $\Gamma_{\bar V}(x)$ where the corresponding  $\Gamma$-epimorphism $\varphi : \Gamma_{\bar V}(x) \rightarrow L$  is one-to-one in restriction to  the elliptic abelian neighbourhood of $H$;} 

$(ii)$ \textit{every proper restricted $\Gamma$-limit group quotient of $\Gamma_{\bar V}(x)$ satisfies the conclusion of the theorem. }

\textit{Then $H$ is $\Gamma$-isomorphic to $\Gamma$. }

\proof  

Since  $\Gamma_{\bar V}(x)$   is a CSA-group either $C_1$ or $C_2$ is malnormal in $H$. We treat the case $C_1$ is malnormal in $H$, the other case being similar.  

Let $D=C_H(C_2)$.  We make the following two assumptions:

$(a)$ $C_1$ and $C_2$ are not conjugate in $H$; 

$(b)$ $D$ is noncyclic;

and we show that we obtain a contradiction. 

By $(a)$ the HNN-extension is seperated, and thus,  since $C_1$ is malnormal, we have $tDt^{-1}=C_{\Gamma_{\bar V}(x)}(C_1)$ and $D= C_{\Gamma_{\bar V}(x)}(C_2)$.  By putting $D'=tDt^{-1}$, and since $\hat H=\<H, \hat C_1, \hat C_2\>=\<H, D'\>$,   we get
$$
\hat H=H*_{C_1}D', \; \Gamma_{\bar V}(x)=\<\hat H, t | D'^t=D\>. 
$$

We notice that $D$ and $D'$ are steal not conjugate in $\hat H$. By construction $D$ and $D'$ are malnormal in $\hat H$. Hence, by Proposition \ref{lem-hnn}, $\Gamma_{\bar V}(x)$ is generated by $\Gamma \cup \{k_1tk_2\}$ for some $k_1, k_2 \in \hat H$.  We replace  $D'$ by $U=D'^{k_1^{-1}}$, $D$ by $V=D^{k_2}$ and $t$ by $r=k_1tk_2$, and thus we get
$$
\Gamma_{\bar V}(x)=\<\hat H, r| U^r=V\>,
$$
and $\Gamma_{\bar V}(x)$ is generated by $\Gamma \cup \{r\}$.  Using normal forms, we conclude that $\hat H=\<\Gamma, U, V\>$. We notice also that $U$ and $V$ are steal not conjugate in $\hat H$.  Using normal forms, we conclude that either $U \cap \Gamma \neq 1$ or $V \cap \Gamma \neq 1$. Without loss of generality we assume that $U \cap  \Gamma \neq 1$.

Clearly $L$ is freely indecomposable with respect to $\Gamma$. If $L$ is $\Gamma$-isomorphic to $\Gamma$ then, since $\varphi$ is one-to-one in restriction to $H$, we get the required conclusion. Hence we assume that $L=\<\Gamma, s|u^s=u\>$.

Since $L$ is generated by $\Gamma \cup\{\varphi(r)\}$, by Proposition \ref{lem-hnn}, $\varphi(r)=\gamma_1s^{\varepsilon}\gamma_2$ for some $\gamma_1, \gamma_2 \in \Gamma$,  $\varepsilon =\pm 1$,  and without loss of generality we assume that $\varepsilon =1$. 

We claim that $U^{\gamma} \leq \<u,s\>$ for some $\gamma \in \Gamma$.  Since $U \cap \Gamma \neq 1$, we let $\gamma_0 \in U \cap \Gamma$. Then $U \leq C_L(\gamma_0)$ and since $U$ is noncyclic, we conclude that $C_L(\gamma_0)=\<u,s\>^{\gamma}$  for some $\gamma \in \Gamma$; and we obtain the required conclusion. 

Replacing $U$ by $U^\gamma$,we assume that $\gamma=1$. Replacing also $V$ by some of its conjugates, we assume that $\gamma_2=1$. Hence we conclude
$$
U=\<u, s^p\>, \; p \in \mathbb Z, \; \; V=U^{\gamma s}. 
$$ 

Suppose towards a contradiction that $p \neq \pm 1$.  We claim that $\hat H=\<\Gamma, U\>*V$, which gives 
$$
\Gamma_{\bar V}(x)=\<\Gamma, U\>*_{U=V^t}\<V^t,t\>=B* \mathbb Z, \; \Gamma \leq B, 
$$ which is a contradiction.

Celarly $\<u, \gamma\>$ is free of rank 2, as otherwise we obtain  that $U$ and $V$ are conjugate in $\hat H$;  which is a contradiction.  It follows that $\<\Gamma, V\>=\Gamma*V$. 

Clearly we also have $\<\Gamma, U\>=\<\Gamma, s^p| u^{s^p}=u\>$.  Since $p \neq \pm 1$, we get  that the length 
$$
|\gamma_1 s^{\varepsilon _1 p}\gamma_2 s^{\varepsilon_2p} \cdots \gamma_n s^{\varepsilon_np}\gamma_{n+1}. s^{-1}\gamma^{-1}d \gamma s|
$$
is greater than 
$$
|\gamma_1 s^{\varepsilon _1 p}\gamma_2 s^{\varepsilon_2p} \cdots \gamma_n s^{\varepsilon_np}\gamma_{n+1}|, \; | s^{-1}\gamma^{-1}d \gamma s|,
$$
for any reduced sequence $\gamma_1 s^{\varepsilon _1 p}\gamma_2 s^{\varepsilon_2p} \cdots \gamma_n s^{\varepsilon_np}\gamma_{n+1}$ of $\<\Gamma, U\>$ and for any nontrivial element $d$ of $U$. Thus by \cite[Lemma 4.2]{AHoucine}, we get the required contradiction. 

Therefore $p=\pm 1$ and thus $U=\<u,s\>$. But this implies $\varphi(H)=L$ and thus $U$ and $V$ are conjugate in $\hat H$; which is also a contradiction.

So finally we conclude that one of the following two cases holds: 

$(1)$ $C_1$ and $C_2$ are conjugate in $H$; 

$(2)$ $C_H(C_1)$ and $C_H(C_2)$ are cyclic. 


We claim that in each case we have 
 $$
 \Gamma_{\bar V}(x)=\<\<\Gamma, c\>, t| c'^t=c\>,
 $$
 where $c' \in \Gamma$ and $H=\<\Gamma, c\>$.

Suppose that $(1)$ holds. By rewritting the HNN-extension, we may assume that $C_1=C_2$. By Proposition \ref{lem-hnn}, $\Gamma_{\bar V}(x)$ is generated by $h_1th_2$ for some $h_1, h_2 \in H$.  Hence 
$$
\Gamma_{\bar V}(x)=\<H, s| D_1^s=D_2\>,
$$
where $s=h_1th_2$, $D_1=C_1^{h_1^{-1}}$, $D_2=C_2^{h_2}$.  Using normal forms we conclude that $H=\<\Gamma, D_1, D_2\>$.

Let $d_1$ (resp. $d_2$) generates $D_1$ (resp. $D_2$).  We claim  that  either $d_1 \in \Gamma$ or $d_2 \in \Gamma$.   Since $d_1$ can be written as a word on $\Gamma \cup \{s\}$, using normal forms we get  $d_1^n \in \Gamma$ or $d_2^n \in \Gamma$ for some $n \in \mathbb Z$, $n\neq 0$. Suppose that $d_1^n \in \Gamma$ for some $n \in \mathbb Z$, $n\neq 0$. Since  $\Gamma$ is existentially closed in $\Gamma_{\bar V}(x)$ (Lemma \ref{lem1}(1)), there exists $\gamma \in \Gamma$ such that $\gamma^n=d_1^n$.   Since $G$ is torsion-free and commutative transitive, we get $d_1=\gamma$ as claimed. Therefore
 $$
 \Gamma_{\bar V}(x)=\<\<\Gamma, c\>, t| c'^t=c\>,
 $$
 where $c' \in \Gamma$ and $H=\<\Gamma, c\>$ as required.
 
 Now suppose that $(2)$ holds with $(1)$ does not hold.  Proceeding as before,  we have 
 $$
 \Gamma_{\bar V}(x)=\<\hat H, r| U^r=V\>,
 $$
where in this case $U$ and $V$ are cyclic and malnormal, and the HNN-extension is seperated. Proceeding as above we conclude that $\hat H=\<\Gamma, U, V\>$ and also that $U \leq \Gamma$ or $V \leq \Gamma$, and without loss of generality we assume that $U \leq \Gamma$.

We claim that $H=\hat H$. Since $\Gamma_{\bar V}(x)$ is a  CSA-group, either $C_1=\<c_1\>$ is malnormal or $C_2=\<c_2\>$ is malnormal in $H$.  We suppose that $C_1$ is malnormal the other case can be treated similalrly. Without loss of generality, after conjugation,  we may also suppose that $C_1\leq U$.  Since $V$ is cyclic we get $c_2=d_1^p$ for some $p \in \mathbb Z$. Hence, proceeding as above, since $c_1$ and $c_2$ are conjugate and $c_1 \in \Gamma$ we conclude that $p=\pm 1$.   Therefore $C_2$ is also malnormal. Finally we conclude that
 $$
 \Gamma_{\bar V}(x)=\<\<\Gamma, c\>, t| c'^t=c\>,
 $$
 where $c' \in \Gamma$ and $H=\<\Gamma, c\>$ as required.

Hence in each case, we have 
 $$
 \Gamma_{\bar V}(x)=\<\<\Gamma, c\>, t| c'^t=c\>,
 $$
 where $c' \in \Gamma$ and $H=\<\Gamma, c\>$.  By Lemma \ref{lem1}(2),  $L=\Gamma_{\bar U}(x)$ for some $U \subseteq \Gamma[x]$. Clearly  $\Gamma_{\bar U}(x)$ is freely indecomposable  relative  to $\Gamma$.  If  $\Gamma_{\bar U}(x)$ is $\Gamma$-isomorphic to $\Gamma$ then we get the required conclusion as $\varphi$ is one-to-one in restriction to $H$.

Therefore   $\Gamma_{\bar U}(x)=\<\Gamma, s|u^s=u\>$.  We claim  that  $\varphi(c) \in \Gamma$ and this will ends the proof.  We steal denote by $c$ the image of $c$ in $\Gamma_{\bar U}(x)$.  Recall that
 $
 \Gamma_{\bar S}(t)=\<\<\Gamma, c\>,t| c'^t=c\>,
 $
 and thus $c'^{\varphi(t)}=c$.  Set $t'=\varphi(t)$. Since $\Gamma \cup\{t'\}$ generates $\Gamma_{\bar U}(x)$, without loss of generality,    $t'=\gamma_1s\gamma_2$ for some $\gamma_1,\gamma_2 \in \Gamma$ by Proposition  \ref{lem-hnn}. 
 
Replacing $s$ by $s^{\gamma_2}$ and $u$ by  $u^{\gamma_2}$ we may assume that $\gamma_2=1$ and we write $t'=\gamma s$.  Therefore $c= s^{-1}\gamma^{-1}c' \gamma s$.  It follows that  $\<\Gamma,c\>=\Gamma*_{u=u^s }\<u^s,  c\>$.

By Claim 1,  $\<u^s,  c\>$ is either free of rank 2 or abelian.  If $\<u^s,  c\>$ is free of rank $2$, then  $\Gamma_{\bar S}(x)$ will be freely decomposable with respect to $\Gamma$; a contradiction to our assumption.  If $\<u^s,  c\>$ is abelian then  $[s,c]=1$ and thus $c= \gamma^{-1}c' \gamma \in \Gamma$ as claimed. This ends the proof of the claim. \qed

Now we treat the two cases of Definition  \ref{weakly-constructible}. 

Case (1).   Let $G=H*_CK$ be the given splitting with $K$ is noncyclic and $\Gamma \leq H$.   Since $G$ is a CSA-group, $C$ is malnormal  either in $H$ or $K$.  We claim that $K$ is abelian. Let 
$$
H'=H*_CC', K'= K, C'=C_K(C),
$$
whenever $C$ is malnormal in $H$ and 
$$
H'=H, K'=C'*_CK, C'=C_H(C),
$$
whenever $C$ is malnormal in $K$. We get $G=H'*_{C'}K'$ with $\Gamma \leq H'$ and $C'<K'$ is malnormal in both $H'$ and $K'$.  By Proposition \ref{lem-hnn},  and without loss of generality, $x=hk$ where $h \in H'$ and $k \in K'$. 

Let $v \in K$. Then $v$ can be written as a reduced word on $\Gamma \cup\{hk\}$. By reducing this word with respect to the structure of the free product with amalgamation, we get  $v \in \<k,C'\>$ and thus $K'=\<k, C'\>$.

By Claim 1,  either $K'$ is abelian or $K'=C'*\<k\>$. Clearly the later case is impossible as otherwise  $\Gamma_{\bar S}(x)$ will be freely decomposable relative  to $\Gamma$; a contradiction with our assumption.  

Therefore $K'$ is abelian and in particular $K$ is abelian as claimed. 
Since $C$ has an infinite index in $K$ we can write $K=C_0 \times C_1 \times \cdots \times C_n$, where $C_0=\<t_0\>$,  $C=\<c\>$ with $c=t_0^p$ for some $p \in \mathbb Z$, $n \geq 1$ and each $C_i$ is cyclic.   We let, for $0 \leq i \leq n$,  
$$
L_0=H*_CC_0, \;L_1=\<L_0, t_1| C_0^{t_1}=C_0\>, \; $$$$L_i=\<L_{i-1}, t_{i}| (C_0\times \dots \times C_{i-1})^{t_i}=(C_0\times \dots \times C_{i-1})\>.
$$

We see that each $L_{i}$ is a proper quotient of $L_{i+1}$ and $\Gamma_{\bar S}(x)=L_n$.  Hence by induction each $L_i$ satisfies conclusions of the theorem for $0 \leq i \leq n-1$. 

We claim  that $L_0$ is $\Gamma$-isomorphic to $\Gamma$. We see that $L_1$ satisfies all the assumptions of Claim 2 and thus $L_0$ is $\Gamma$-isomorphic to $\Gamma$ as desired. 
 
We claim that $n=1$. Suppose towards a contradiction that $n \geq 2$.

 We have  
$$
L_2=\<\<L_0, t_1| t_0^{t_1}=t_0\>, t_2| t_0^{t_2}=t_0, t_1^{t_2}=t_1\>.  
$$

By Proposition \ref{lem-hnn},  $L_2$ is generated by $\Gamma \cup \{h_1t_2h_2\}$ where $h_1,h_2 \in L_1$. Again, since $L_1$ is generated by $\Gamma \cup\{h_1h_2\}$ we find, by Proposition \ref{lem-hnn}, $h_1h_2= \gamma_1t_1^{\pm 1}\gamma_2$ for some $\gamma_1, \gamma_2 \in L_0$. Now there exists a word $w(\bar x;y)$ such that $t_1=w(\bar \gamma; (h_1t_2h_2))$,  and thus in $L_0 \times \<t_1\>\times \<t_2\>$ we have 
$$
w(\bar \gamma; h_1t_2 h_2)=w(\bar \gamma; h_1h_2t_2)=v(\bar \gamma; \gamma_1 \gamma_2)(t_1^{\pm 1}t_2)^p=t_1,
$$
for some $p\neq 0$, $p\in \mathbb Z$, which is clearly a contradiction.

Hence $n=1$ as claimed and finally   
$$
\Gamma_{\bar s}(x)=\<\Gamma,  t_1| t_0^{t_1}=t_0\>. 
$$

Case $(2)$.  Let $\Gamma_{\bar S}(x)=\<H,t |c_1^t=c_2\>$.  $\Gamma_{\bar S}(x)$ satisfies all the assumptions of Claim 2, and thus $H$ is $\Gamma$-isomorphic to $\Gamma$.   Therefore $\Gamma_{\bar S}(x)=\<\Gamma, t| c'^t=c\>$.  Since  $\Gamma$ is existentially closed in $\Gamma_{\bar S}(x)$,   $c$ and $c'$ are conjugate in $\Gamma$. Thus $\Gamma_{\bar S}(x)$ can be rewritten as $\<\Gamma, s| u^s=u\>$ and we obtain the required conclusion.  This ends the proof in this case and the proof of the theorem. \qed

\smallskip
\noindent
\textbf{Proof of Theorem \ref{thmmain2}.}  

Theorem \ref{thmmain} shows $(2) \Rightarrow (1)$, so we show   $(1) \Rightarrow (2)$.  Let $H=\<a,b\>$ be a nontrivial two-generator subgroup of $\Gamma$.  We may suppose without loss of generality that $a$ is root-free.  We claim that the group $\Gamma*_a\<a, b'\>$ is a restricted $\Gamma$-limit group, where $\<a, b\> \cong \<a,b'\>$.  Since $a$ is root-free, by applying \cite[Lemma 5.4]{groves-2007}, we see that the group $\<\Gamma, t | a^t=a\>$ is a restricted $\Gamma$-limit group. We have  $\<\Gamma, b^t\>=\Gamma*_a\<a,b^t\>$ with $\<a,b^t\> \cong \<a,b\>$.  Hence $\Gamma*_a\<a,b'\>$ is a restricted $\Gamma$-limit group, where we can take $b'=b^t$.  

Let
$$
S(x)=\{w(x) \in \Gamma[x]| \Gamma*_a\<a,b'\> \models w(b')=1 \}. 
$$

It is not hard to see that 
$$
\Gamma_{\bar S}(x) \cong_{\Gamma} \Gamma*_a\<a,x\>,
$$
with $\<a,x\> \cong \<a,b\>$. 

Suppose that $[a,b] \neq 1$. Then  $V(S)$ is infinite and irreducible.  Hence, by (1), $V(S)$ is a coset of a centralizer. So let $u,v \in \Gamma$ such that $V(S)=v C_\Gamma(u)$.  By applying \cite[Lemma 5.4]{groves-2007}, we conclude that
$$
\Gamma_{\bar S}(x) \cong_{\Gamma} \<\Gamma, s | u^s=u\>,
$$
where $x=vs$. 

Suppose that $\<a,x\>$ is not free of rank $2$. Then  there exists a nontrivial relation and using normal forms, we conclude that either $a^p=u^q$ or $v^{-1}a^pv=u^q$ for some $p, q \in \mathbb Z$.  If the latter case holds then we may replace $u$ by $vuv^{-1}$ and $s$ by $vsv^{-1}$ and thus we get $x=sv$. Thus we conclude that we may assume $a^p=u^q$  for some $p, q \in \mathbb Z$ and $x=vs$ or $x=sv$. Since $a$ and  $u$ are root-free, we get $a=u^{\pm 1}$ and without loss of generality, we assume that $a=u$. 

Returning to our first construction, we get 
$$
\<\Gamma, b^t\>=\Gamma*_a\<a,b^t\>=\<\Gamma, s|a^s=a\>\leq  \<\Gamma, t|a^t=a\>, 
$$
and thus $s=a^pt^q$  for some $p, q \in \mathbb Z$. Hence $b^t= v a^pt^q$ or $b^t=t^q a^p v$. In the group $\Gamma \times \<t|\>$ we get   $q=0$.  Hence we find $b^t=va^p$ and thus $b \in \<a\>$,  which is a contradiction. Therefore $\<a,x\>$ is free of rank $2$. \qed

\section{Appendix}

In this appendix,  we give a proof of the following theorem, where $\Gamma$ is steal a torsion-free hyperbolic group. For the notions used here, and which are not defined, we  refer  the reader to \cite{Champ-Guirardel}. 

\begin{thm} \label{athm1}Let $G$ be a restricted $\Gamma$-limit group. If $G$ is not $\Gamma$-isomorphic to $\Gamma$ and if it is freely indecomposable relative to $\Gamma$,  $G$   is weakly constructible. \qed
\end{thm}

\begin{defn}  A cyclic splitting (relative to $\Gamma$) is \textit{essential} if any edge group is of infinite index in any vertex group. \qed
\end{defn}

The proof of the following proposition is similar to that of \cite[Theorem 3.7]{groves-2007} and it is left to  the raider. 

\begin{prop} \label{aprop1}A restricted $\Gamma$-limit group which  is not $\Gamma$-isomorphic to $\Gamma$ and which is freely indecomposable relative to $\Gamma$ admits an essential cyclic splitting (relative to $\Gamma$). \qed
\end{prop}

\begin{defn} \cite[Definition 8.3]{Sela-Diophan1} Let $G$ be  a restricted $\Gamma$-limit group  which  is not $\Gamma$-isomorphic to $\Gamma$ and which is freely indecomposable relative to $\Gamma$. The \textit{restricted modular group $RMod(G)$} is the subgroup of $Aut(G)$ generated by the following families of automorphisms of $G$, which fixe pointwise  the vertex group stabilized by $\Gamma$ in the restricted cyclic JSJ-decomposition of $G$ with respect to $\Gamma$: 

$(1)$ Dehn twists along edges of the restricted cyclic JSJ-decomposition of $G$. 

$(2)$ Dehn twists along essential s.c.c. in CMQ vertex groups in the restricted cyclic JSJ-decomposition of $G$.

$(3)$ Let $A$ be an abelian vertex group in the restricted cyclic JSJ-decomposition of $G$. Let $A_1<A$ be the subgroup generated by all edge groups connected to the vertex stabilized by $A$ in the cyclic JSJ-decomposition of $G$. Every automorphism of $A$ which fixes pointwise $A_1$ can be extended to an automorphism of $G$ which fixes the vertex stabilized by $\Gamma$. We call these generalized Dehn twists and they form the 
third family of automorphisms that generate $RMod(G)$. \qed
\end{defn}

\begin{defn} (Shortening quotients)  Let $G$ be a restricted $\Gamma$-limit group endowed with a finite generating set $B$.

$(1)$  A $\Gamma$-homomorphism $h : G \rightarrow \Gamma$ is said \textit{short}   if  
$$
\max_{b \in B}|h(b)|\leq \max_{b \in B}|h(\tau(b))|,
$$
for any restricted modular automorphism $\tau \in RMod(G)$. Here $|.|$ denotes  the word length with respect to some fixed,  for all the rest of this section, finite generating set of $\Gamma$. 

$(2)$ Let $(h_n : G  \rightarrow \Gamma)_{n \in \mathbb N}$ be a sequence of short $\Gamma$-homomorphisms. The group $G/ Ker_{\infty}((h_n)_{n \in \mathbb N})$ is called a \textit{shortening quotient} of $G$. \qed
\end{defn}

\begin{thm}\label{athm2} \cite[Claim 5.3]{Sela-Diophan1}\cite[Proposition 1.15]{Sela-hyp} Let $G$ be  a restricted $\Gamma$-limit group  which  is not $\Gamma$-isomorphic to $\Gamma$ and which is freely indecomposable relative to $\Gamma$. Then every shortening quotient of $G$ is a strict quotient. \qed
\end{thm}

\begin{prop} \label{aprop2} Let $G$ be a restricted $\Gamma$-limit group  which  is not $\Gamma$-isomorphic to $\Gamma$ and which is freely indecomposable relative to $\Gamma$. Then either $G$ is a free extension of a centralizer or $G$ is weakly constructible. 
\end{prop}

We begin first with the following lemma which is analogous to \cite[Proposition 4.12]{Champ-Guirardel}. 

\begin{lem} \label{alem1}Let $H$ be a restricted $\Gamma$-limit group with a one edge cyclic splitting $H=A*_CB$ or $A*_C$ satisfying the following property:  there exists a $\Gamma$-epimorphism $\varphi : H \rightarrow L$, where $L$ is a restricted $\Gamma$-limit group, such that $\varphi$ is one-to-one in restriction to  the elliptic abelian neighbourhood of each  vertex group. 

Then there exists a sequence of Dehn twists $(\tau_i)_{i \in \mathbb N}$ on $H$, fixing pointwise $A$,  such that $(\varphi \circ \tau_i)_{i \in \mathbb N}$ converges to the identity of $H$. 
\end{lem}

\proof  Proceeding as in the proof of Theorem \ref{thm2}, one first transform the given splitting to another one which is either 1-acylindrical or a free extension of a centralizer. Then the rest of the  proof proceeds in a similar way to that of \cite[Proposition 4.12]{Champ-Guirardel} and \cite[Theorem 5.12]{Sela-Diophan1}, by using \cite[Lemma 5.4]{groves-2007} instead of Baumslag's lemma \cite[Lemma 3.5]{Champ-Guirardel} and by choosing the Dehn twist along some $c \in C$.  \qed

\smallskip
\noindent\textbf{Proof of Proposition \ref{aprop2}.} 

The proof proceeds in a similar way to that of \cite[Proposition 4.18, Proposition 4.18]{Champ-Guirardel} and \cite[Proposition 5.10]{Sela-Diophan1}.  We suppose that $G$ is not a free extension of a centralizer and that it does not satisfy (1) of Definition \ref{weakly-constructible}.

Let $\Lambda$ to be the cyclic JSJ-decomposition of $G$ which is nontrivial by Proposition \ref{aprop1}. Suppose first that $\Lambda$ has  an abelian vertex group  $G_v$ such that $A_1$ has an infinite index in $G_v$, where $A_1$ is the group generated by incident edge groups. Then in that case $G$ can be written as a nontrivial free extension of a centralizer.

Thus we may assume that for each abelian vertex group $G_v$, the subgroup generated by incident edge groups has finite index in $G_v$.  Hence by definition each restricted modular automorphism $\tau$ is a conjugation in restriction to each nonsurface type vertex group,  to each edge group and the identity on the vertex group containing $\Gamma$. 

Let $(f_i : G \rightarrow \Gamma)_{i \in \mathbb N}$  be a sequence of $\Gamma$-homomorphisms  converging to the identity of $G$.  For each $i \in \mathbb N$ choose $\tau_i$ to be a restricted modular automorphism such that $f_i \circ \tau_i$ is short.  Up to extracting a subesequence,  we may assume that $(f_i \circ \tau_i)_{i \in \mathbb N}$ converges to a restricted $\Gamma$-limit group $L$ and we let $\varphi :G \rightarrow L$ to be the natural map.  By Theorem \ref{athm2}, $L$ is a proper quotient.

Proceeding as in \cite[Proposition 4.18]{Champ-Guirardel}, we conclude that $\varphi$ is one-to-one in restriction  to the elliptic abelian neighboorhood of each  nonsurface vertex group and of the vertex group containing $\Gamma$

Let $A$ be the vertex group containing $\Gamma$ and  let $e$ be an edge incident to $A$. Write $H=A*_CB$ or $H=A*_C$ the subgroup of $G$ corresponding to the amalgam or HNN-extension carried by $e$.

Suppose that $H=A*_CB$ and $B$ is abelian. Since  $H \leq \hat A$,  it follows that  $\varphi$ is one-to-one in restriction to $H$. Let $\bar \Lambda$ be the graph of groups obtained by collapsing $e$. Then $H$ is a vertex group and $\varphi$ is one-to-one in restriction to elliptic abelian  
neighboorhood of each vertex group of $\bar \Lambda$. If there is a another vertex abelian group $H'$ connected to $H$, we do the same construction. At the end of the procedure we get a cyclic splitting $\Lambda'$ such  that if an  edge is connected to the vertex group containing $\Gamma$ in $\Lambda'$ with different end points then in the corresponding amalgam  $H=A*_CK$, $K$ is nonabelian. But this contradicts our hypothesis; because in that case $G=D_1*_CD_2$ whith $D_2$ noncyclic and $\Gamma \leq D_1$. 

Hence  $G$ can be written as   $G=\<K, t_1, \dots, t_n | C_i^{t_i}=C'_i\>$  for some cyclic subgroups   $C_1, \dots, C_n, C'_1, \dots, C'_n$ of $K$ and $\Gamma  \leq A \leq K \leq \hat A$ and $\varphi$ is one-to-one in restriction to $\hat A$. If for some $i$ and $a \in  K$, $C_i \cap {C'_i}^a\neq 1$  then $G$ can be written as a free extension of a centralizer. Hence for any $i$ and $a \in K$, $C_i \cap {C'_i}^a= 1$.  

If $n=1$ we get the required conclusion. So we suppose that $n \geq 2$. 

Let $H=\<K, t_1|C_1^t=C'_1\>$ and let  $e$ be the edge  corresponding to $C_1$. Lemma \ref{alem1} applies in this case and we get a sequence of  Dehn twists $(\tau_i)_{i \in \mathbb N}$ on $H$ such that $(\varphi \circ \tau_i)_{|H}$ converges to the identity of $H$.  Up to exctracting a subsequence, we may assume that $(\varphi \circ \tau_i)_{i \in \mathbb N}$ converges  to a $\Gamma$-epimorphism $\phi : G \rightarrow L'$, where $L'$ is a restricted $\Gamma$-limit group and where we identify $\tau_i$ with its natural extension to the entire group $G$. Let $\bar \Lambda$ be the graph of groups obtained by collapsing $e$.  By construction $\phi$ is one-to-one in restriction to the elliptic abelian neighboorhood of the vertex group. If the obtained $\phi$ is not one-to-one, we conclude by induction on $n$. 


So suppose that $\phi$ is one-to-one. We consider in this case the  connected component $\Lambda_1$ of $\Lambda\setminus e$. Then $\varphi$ is one-to-one in restriction to the elliptic abelian neighboorhood of the fundamental group of $\Lambda_1$. Hence we obtain a one edge cyclic splitting of $G$ such that $\varphi$ is one-to-one in restriction to the elliptic abelian neighboorhood of the vertex group. \qed


\noindent\textbf{Proof of Theorem  \ref{athm1}.}

Let $G$ be a restricted $\Gamma$-limit group which is not $\Gamma$-isomorphic to $\Gamma$ and which is freely indecomposable relative to $\Gamma$. By the descending chain condition on restricted $\Gamma$-limit groups, we may asssume that every  restricted $\Gamma$-limit proper quotient of $G$ satisfies the conclusion of the theorem if it satisfies its hypothesis.

By Proposition \ref{aprop2}, we may assume that   $G$ is a nontrivial free extension of a centralizer. Set $G=\<H,t|C^t=C\>$ where $C$ is a nontrivial abelian subgroup of $H$ and $\Gamma \leq H$.  

Define $\phi : G \rightarrow H$ by $\phi(t)=1$ and the identity on $H$.  

If $C$ is cyclic then  we get the required conclusion. So we suppose that $C$ is noncyclic.  Clearly $H$ is not $\Gamma$-isomorphic to $\Gamma$.  Similarly if $H$ is freely decomposable with respect to $\Gamma$ then  $C$ is contained in some conjugate of a factor and thus $G$  is itself freely decomposable with respect to $\Gamma$.

Hence $H$ satisfies the hypothesis of the theorem and by induction we conclude that $H$ is weakly constructible. 

Suppose that $H=A*_TB$ where $\Gamma \leq A$, $T$ is nontrivial and cyclic  and $B$ is noncyclic.  Since $C$ is noncyclic we conclude, up conjugation,  that $C \leq A$ or $C \leq B$.  Therefore $G$ can be written as $A'*_TB'$ with $B'$ is noncyclic.

 Now suppose that $H=A*_T$. Let $L$ be the proper restricted $\Gamma$-limit  quotient of $H$ given by the definition and  let $\varphi : H \rightarrow L$ be the corresponding $\Gamma$-epimorphism. 
 
 Suppose first that $C$ is not elliptic in the splitting  $H=A*_T$. Since $C$ is noncyclic, we conclude that $H$ can be written $H=A*_TC'$ where $C'$ is a conjugate of $C$. Hence $G$ can be written as $G=A*_TC''$, with $C''$ is noncyclic and we get the required conclusion.  
 
Suppose now  that  $C$ is elliptic in the splitting $H=A*_T$ and without loss of generality that $C \leq A$.

Let $C'=C_L(\varphi(C))$ and let $L'= \<L, s|C'^s=C'\>$.  Then $L'$ is a restricted $\Gamma$-limit group. Define $\varphi' : G \rightarrow L'$ by $\varphi'_{|H}$ to be $\varphi$ and $\varphi'(t)=s$.  Now $L'$ is a strict quotient of $G$ as $L$ is a proper quotient of $H$.

Then $G=\<A,t\>*_T=\<A,t|C^t=C\>*_T$ with $T \leq A$. Hence $G$ has a cyclic splitting and with $\varphi'$ is one-to-one  in restriction to  the elliptic abelian neighboorhood of $\<A,t|C^t=C\>$. \qed

\smallskip
We close this appendix with  the following proposition. 

\begin{prop}  Let $\Gamma$ be a torsion-free hyperbolic group and let  $H \leq \Gamma$ be a proper subgroup definable by a quantifier-free formula. Then $H$ is abelian. 
\end{prop}

\proof 

Since $\Gamma$ is equationally noetherian, $H$ is closed in the Zariski topology. Hence $H$ is definable by a finite union of varieties. Without loss of generality we assume that $H$ is definable by an equation, the general case can be treated similarly. So suppose that $H$ is definable by $w(c_1, \dots, c_p;x)=1$.  Let $a \in H$.  By \cite[Lemma 5.4]{groves-2007}, since for any $n \in \mathbb N$,  $w(c_1, \dots, c_p; a^n)=1$,  we obtain  $a \in C_\Gamma(c_1) \cup \dots \cup C_{\Gamma}(c_p)$.  But if $H$ is nonabelian,  $H$ contains a nonabelian  free subgroup and we get a contradiction. \qed


\bibliographystyle{alpha}
\bibliography{biblio}
\end{document}